\input amstex
\input xy
\xyoption{all}
\documentstyle{amsppt}
\document
\magnification=1200
\NoBlackBoxes
\nologo
\hoffset1.5cm
\voffset2cm
\pageheight {16cm}
\def\A{\Cal{A}}

\def\M{\Cal{M}}
\def\Q{\Cal{Q}}
\def\O{\Cal{O}}
\def\r{\roman}
%\hfill{\it file Amywork/Gr-Ver-operads.tex, version 22.10.2016}

\bigskip

\centerline{\bf GROTHENDIECK--VERDIER DUALITY PATTERNS}

\medskip

\centerline{\bf IN QUANTUM ALGEBRA}

\bigskip

\centerline{\bf Yuri I.~Manin}

%\medskip

%\centerline{Max--Planck--Institut f\"ur Mathematik, Bonn, Germany}

\bigskip

{\bf Abstract.} After a brief survey of the basic definitions of the Grothendieck--Verdier
categories and dualities,  I consider in this context introduced earlier dualities in the categories
of quadratic algebras and operads, largely motivated by the theory of quantum groups.
Finally, I argue that Dubrovin's ``almost duality'' in the theory of Frobenius
manifolds and quantum cohomology also must fit a (possibly extended)
version of Grothendieck--Verdier duality.

\bigskip


\medskip{


\bigskip

\centerline{\bf Introduction and summary}

\medskip

Duality is one of the omnipresent and elusive kinds of symmetries
in mathematics. Albert Einstein allegedly joked that
his most important mathematical discovery was summation over
repeating sub-/superscripts in tensor analysis. It was a wise joke.
\smallskip

The standard linear duality functor establishes an equivalence
between the category of finite--dimensional linear spaces
(over a fixed field $k$) and its opposite category.
It is replaced by Serre's duality, if one passes to the category
of vector bundles over a smooth projective manifold; and
by Grothendieck--Verdier duality if one passes to the derived
or enhanced triangulated categories of sheaves in a more
general context.

\smallskip

Here I am interested in the interaction of duality with monoidal 
structure(s), such as appearance of black and white products
related by the Koszul duality
in the category of quadratic algebras ([16]).
M.~Boyarchenko and V.~Drinfeld developed  the general duality formalism
in non necessarily symmetric monoidal categories in [4], and 
I want to look at several of the constructions suggested earlier
in the light of their formalism.

\smallskip

Sec.~1 is a brief survey of the relevant part of [4] supplied by the discussion of one of their
results in a slightly more general setting (see comments after Proposition 1.2).

\smallskip

In Sec.~2, I show that  old constructions of quantum groups via
category of quadratic (and in fact, more general) algebras in [16], and quadratic  operadic
duality in [11],  fit in the (extended) context
of [4].

\smallskip

Finally in Sec.~3, I argue that constructions of Dubrovin's ``almost duality'' [10] in the theory of Frobenius 
manifolds and quantum cohomology suggest their categorical enhancement
that might also fit in a Grothendieck--Verdier context.

\smallskip

{\it Acknowledgments.} I am grateful to Vladimir Drinfeld and Bruno Vallette whose works
and comments on earlier drafts of this note were very useful for me.

\smallskip

Finally, I am very happy to dedicate this article to the anniversary of the Izvestiya RAN,
where sixty years ago  my first paper  was published.

\bigskip

\centerline{\bf 1. Dualities in monoidal categories}

\medskip

{\bf 1.1. GV--categories.} Consider a monoidal category $(\M, \bullet , \bold{1})$. (For brevity, we omit
here notation for the relevant associativity and unit morphisms and constraints, cf. e.~g. [15], p.~580). 
\smallskip
The basic definitions in [4] introduce
two notions:

\smallskip

(i) An antiequivalence functor  $D_K:\, \M\to \M^{op}$ is called {\it duality functor}
if there exists such an object $K$ of $\M$ that
 for every object $Y$ in $\M$, the functor 
$X\mapsto \roman{Hom} (X\bullet Y,K)$ is representable by 
$D_KY$. In this case
$K$  is called {\it dualizing object} for $\M$.

\smallskip

(ii) The datum $(\M,\bullet , \bold{1},K)$ as above is called {\it a Grothendieck--Verdier category,} GV--category for short.
   
 \medskip
 
 As was noted in [4], GV--categories were also studied in the literature
 under the name {\it *--autonomous categories}, cf. [1], [2].  
 Clearly, researcher groups studying this kind of duality, were motivated
 by different classes of examples and arguments: say, derived categories of schemes
 and stacks for [4], and normed spaces for [1]. In [11], a duality functor $\bold{D}$ was constructed
 also on the category of $dg$--operads and explicitly compared to the Verdier
 duality ([11], p.~205).

 \smallskip

 For us, one of the most interesting properties of GV--categories is the following one
 ([4], sec. 4.1). 
 
 \smallskip

{\bf 1.2. Proposition.} {\it  (i) Let  $(\M, \bullet , \bold{1}, K)$ be a GV--category.  Then, using $D_K$, one can define a new monoidal structure
 $\circ : \M\times \M\to \M$  putting
 $$
 X\circ Y:= D_K^{-1} (D_KY\bullet D_KX).
 \eqno(1.1)
 $$
 Its unit object is $K$.
 
 \smallskip
 
 (ii)  If $\M$ is an r--category that is,
 if $\bold{1}$ is a dualizing object, then 
 the monoidal products $(\bullet , \circ )$ are connected by the canonical functorial morphisms
 $$
 X\bullet Y\to X\circ Y
\eqno(1.2)
 $$
 compatible with the respective associativity constraints}
  
 \medskip
 
{\it Comments.} (a) Notice that our $(\bullet, \circ )$ correspond to $(\otimes ,\odot )$ in [4] respectively.
 
 \smallskip
 
 (b) Examples of quadratic algebras ([16]) and binary quadratic operads ([11]) show that
 Proposition 1.2 holds for their respective categories, although $\bold{1}$ in neither
 of two cases is a dualizing object (cf. remark at the end of subsection 2.2 below). 
For this reason I will briefly repeat here those
 arguments of [4] that do not use the latter requirement.

 \smallskip
 
  (c) Let us produce first the identity
  isomorphisms in $(\M, \circ )$. In [4], (2.9), it is proved that
  $D_KK$ is canonically isomorphic to the unit object $\bold{1}$ for $\bullet$. Therefore
  $$
  K\circ X= D_K^{-1}(D_KX\bullet D_KK)=D_K^{-1}(D_KX)=X,
  $$
  $$
  X\circ K=D_K^{-1}(D_KK\bullet D_KX)= D_K^{-1}(D_KX)=X.
  $$
   
  \smallskip
  (d) Assuming the existence of (1.2), let us define the associativity 
  isomorphisms in  $(\M, \circ ,K)$.
  
    \smallskip
  Start with associativity morphisms in  $(\M, \bullet , \bold{1})$: for all objects $X,Y,Z$ we
  have isomorphisms
  $$
  \alpha (X,Y,Z): \  (X\bullet Y)\bullet Z \to X\bullet (Y\bullet Z).
  \eqno(1.3)
  $$
  Using the classical Mac Lane theorem,  we may and will assume that  $(\M, \bullet, \bold{1})$  
  is {\it strict} monodical category, and $D_K$ is an isomorphism. Then all associativity (and identity) morphisms 
  for $\bullet$ are identities.
  Using this and  (1.1), we can identify the first objects of the following two lines:
  $$
  D_K(X\circ (Y\circ Z))=D_K(Y\circ Z)\bullet D_KX = (D_KZ\bullet D_KY)\bullet D_KX,
  $$
  $$
  D_K((X\circ Y)\circ Z)=D_KZ\bullet D_K(X\circ Y)=D_KZ\bullet (D_KY\bullet D_KX).
  $$
  Applying finally $D_K^{-1}$, we get the associativity morphisms in $(M,\circ , K)$
  $$
  \beta (X,Y,Z) :\   (X\circ Y)\circ Z \to  X\circ (Y\circ Z).
  $$

  \bigskip
  
  \centerline{\bf 2. GV--categories of quadratic algebras}
  
  \smallskip
  
\centerline{\bf and binary quadratic operads}
  
  \medskip
  
  {\bf 2.1. Quadratic algebras.} The category of quadratic algebras $\Q\A$
  over a field $k$ was defined in [16] in the following way.
  
  \smallskip
  
  {\it An object} $A$ of $\Q\A$ is an associative graded algebra $A=\oplus^{\infty}_{i=0}A_i$
  such that $A_0=k$, $A_1$ generates $A$ over $k$; and finally, the ideal of all relations
  between elements in $A_1$ is generated by its subspace of quadratic relations
  $R(A)\subset A_1^{\otimes 2}.$
  
  \smallskip
  
  As in [16], I will briefly summarise this description as $A \leftrightarrow \{A_1, R(A)\}$.
  
  \smallskip
  
  {\it A morphism} $f:\,A\to B$ is simply a morphism of graded algebras identical on $k$.
  Clearly, morphisms $A\to B$ are in bijection with the set of linear maps $f_1:\,A_1\to B_1$
  such that $(f_1\otimes f_1)(R(A))\subset R(B).$
  
  \smallskip
  
  Put 
  $$
  A\bullet B \leftrightarrow  \{ A_1\otimes B_1, S_{23} (R(A)\otimes R(B))\}
  \eqno(2.1)
  $$
  where 
  $$
  S_{23}:\,  A_1^{\otimes 2}\otimes B_1^{\otimes 2}\to (A_1\otimes B_1)^{\otimes 2} 
  \eqno(2.2)
  $$ 
  sends $a_1\otimes a_2\otimes b_1\otimes b_2$ to $a_1\otimes b_1\otimes a_2\otimes b_2$.
  
  \smallskip
  
  Finally, given a quadratic algebra $A$, define its dual algebra $A^!$ by the convention
  $$
  A^! \leftrightarrow \{A_1^*, R(A)^{\perp}\}  .
  \eqno(2.3)
  $$
  Here $A_1^*$ is the space of linear functions $A_1\to k$, and $R(A)^{\perp}$ is
  the subspace orthogonal to $R(A)$ in $(A_1^*)^{\otimes 2}$.
  \medskip
  
  {\bf  2.2. Proposition.} {\it (i) $(\Q\A ,\bullet ,\bold{1}:=k[\tau ]/(\tau^2))$ is the monoidal category
  whose unit is the  algebra of ``dual numbers''.
  \smallskip
  
  (ii) The map $A\mapsto A^{!}$ extends to the duality functor
  $D_K:\, \Q\A \to \Q\A^{op}$ with dualizing object
  $$
  K:=k[t] = \bold{1}^!
  \eqno(2.4)
  $$
  so that $(\Q\A,\bullet , \bold{1}, K)$ is a GV--category.
  
  \smallskip
  
  (iii) The respective white product (1.1) is given explicitly by
  $$
  A\circ B \leftrightarrow \{A_1\otimes B_1, S_{23}(R(A)\otimes B_1^{\otimes 2} +A_1^{\otimes 2}\otimes R(B))\} .
  \eqno(2.5)
  $$
  Its unit object is K,
  whereas the morphism (1.2) is induced by the obvious embedding
  $$
  R(A\bullet B)= S_{23} (R(A)\otimes R(B)) \subset S_{23}(R(A)\otimes B_1^2 +A_1^{\otimes 2}\otimes R(B)) = R(A\circ B ).
  $$
 }
 \smallskip
 
 {\bf Sketch of proof.} This Proposition is just a reformulation in the language of [4]
 of the results stated and proved in [16] , pp. 19--28.
 
 \smallskip
 
 The key requirement in the definition of duality functor by Boyarchenko and Drinfeld is the functorial
 isomorphism $\r{Hom} (X\bullet Y,K)\cong \r{Hom} (X, D_KY)$. 
 \smallskip
  In $\Q\A$, it follows from Theorem 2 on p. 25 of [16] whose relevant special case gives
  $$
 \r{Hom}_{\Q\A}(A\bullet B, K)\cong \r{Hom}_{\Q\A}(A, B^!\circ K) \cong
 \r{Hom}_{\Q\A}(A, B^!)
 $$
 because $K$ is the unit object for $(\Q\A ,\circ )$. Thus, $D_K=!$,
 and  $D_K$ is an equivalence $\A\Q\to \A\Q^{op}$, because $!\cdot ! \cong \r{Id}_{\Q\A}$.
 \smallskip
 
 We leave the remaining details to the reader.
 
 \medskip
 {\it Remark.} Here I will show that $\bold{1}$ is not a dualizing object in $\Q\A$.
 In fact,
 $$
 \r{Hom}_{\Q\A}(A\bullet B, \bold{1})\cong \r{Hom}_{\Q\A}(A, B^!\circ \bold{1}) .
 $$
 But the functor $\Q\A\to \Q\A^{op}$ acting on objects as $B\mapsto B^!\circ \bold{1}$
 is not an equivalence, because from (2.3) and (2.5) it follows that
 $$
 B^!\circ \bold{1} \leftrightarrow \{B_1^*, (B_1^*)^{\otimes 2}\}.
 $$
  \medskip
  
   {\bf 2.3. Binary quadratic operads.} Extending [16], V.~Ginzburg and M.~Kapranov ([11]) defined the category of binary quadratic operads $\Q\O$ with black and white products  in the following way.
   
   \smallskip
   
   As above, $k$ is a ground field.  Arity--components of objects $Op \in \Q\O$ are finite--dimensional
   linear spaces over $k$ such that $Op(1)=k$. Moreover, $Op$ must be generated by the space
   $E:=Op(2)$ with structure involution  $\sigma$. The sign $\otimes$ in this subsection will always denote
   the tensor product over $k$.
    \smallskip
   
   Generally, linear spaces of generators of associative algebras are replaced in various categories of operads
   by collections of $S_n$--modules. In the category $\Q\O$, the following collections are considered:
$ \{ E(2):=E, E(n)=0, n>2\}$. Such collections will play the role
   of the linear spaces of generators $A_1$ in the formalism of quadratic algebras.
   Free associative algebra generated by $A_1$ is replaced now by the free operad,
   generated by the respective collection.  In our case it will be denoted $F(E)$
   (see its explicit description in [GiKa94], (2.1.8)).  Here $E$ can be an arbitrary $S_2$--module.
   
   \smallskip
   
   All relations between generators must be generated by 
   an $S_3$--invariant subspace $R\subset F(E)(3)$.
   
 \smallskip
  
  As in [16], I will briefly summarise this description as 
  $$
  Op \leftrightarrow \{E, R\}
  \eqno(2.6)
  $$
  \smallskip
 Operadic morphisms are defined as in [11],  (1.3.1).
 
 \smallskip
 
 Now define white and black products in $\Q\O$ by the following data.
 For $j=1,2$, let $Op_j\leftrightarrow (E_j, R_j)$. 
 
 \smallskip
 
 As is explained in [11, Erratum], for all $n\ge 2$ we have canonical maps
 $$
 \varphi_n:\ F(E_1\otimes E_2)(n)\to F(E_1)(n)\otimes F(E_2)(n),
 $$
 $$
 \psi_n:\  F(E_1)(n)\otimes F(E_2)(n) \to  F(E_1\otimes E_2)(n).
 $$
 Then, by definition, 
$$
  Op_1\bullet Op_2 \leftrightarrow  \{ E_1\otimes E_2, \psi_3 ((R_1\otimes F(E_2)(3))\cap (F(E_1)(3)\otimes R_2)\}
  \eqno(2.7)
  $$
  This is an analog of the formula (2.1). However, the appearance of the maps $\varphi_n$ and $\psi_n$
  is an essential new element in the operatic framework.

      \smallskip
 Similarly, the white product is defined by
 $$
  Op_1\circ Op_2 \leftrightarrow  \{ E_1\otimes E_2, \varphi_3^{-1} ((R_1\otimes F(E_2)(3)) + (F(E_1)(3)\otimes R_2)\}
  \eqno(2.8)
  $$

    \smallskip
  
  Finally, given a quadratic operad (2.6), define its dual operas $Op^!$ by the convention
  $$
  Op^! \leftrightarrow \{E^*, R^{\perp}\}  .
  \eqno(2.9)
  $$
  Here $E^*$ is the space of linear functions $E\to k$, and $R^{\perp}$ is
  the subspace orthogonal to $R$ in $F(E^*)(3)$. An important additional
  remark is this: since $E$ is the component of arity 2, it is endowed
  with an action $\rho$ of $S_2$. Respectively, $E^*$ is endowed with
  the action $\rho^*\otimes \r{Sgn}$ where $Sgn$ is the sign--representation.
  \smallskip

  We have now the following analog of the Proposition 2.2. Denote by $Lie$,
  resp. $Comm$ the operads classifying Lie, resp. commutative $k$--algebras (possibly, without unit)
  
  \medskip
  
 {\bf  2.4. Proposition.} {\it (i) $(\Q\O ,\bullet ,\bold{1}:=Lie)$ is a symmetric monoidal category
  whose unit is the  operad $Lie$.
  \smallskip
  
  (ii) The map $Op\mapsto Op^{!}$ extends to the duality functor
  $D_K:\, \Q\O \to \Q\O^{op}$ with dualizing object
  $$
  K:= Comm  = \bold{1}^!
  \eqno(2.10)
  $$
  so that $(\Q\O,\bullet , \bold{1}, K)$ is a GV--category.
  
  \smallskip
  
  (iii) The respective white product (1.1) is given explicitly by (2.8).
\smallskip
  Its unit object is K,
  whereas the morphism (1.2) is induced by the  embedding
  $$
  \psi_3 ((R_1\otimes F(E_2)(3)) \subset
\varphi_3^{-1} ((R_1\otimes F(E_2)(3)) + (F(E_1)(3)\otimes R_2)
  $$
  
 }
 \medskip
 
 The proof uses the same arguments as the proof of Proposition 2.2.
 The key adjointness identity 
 $$
 \r{Hom}_{\Q\O}(Op_1\bullet Op_2, Op_3)=  \r{Hom}_{\Q\O}(Op_1, Op_2^!\circ Op_3)
 $$
 and identification of unit and dualizing objects are sketched in [GiKa94].

\medskip

 {\bf 2.5.  2--monoidal categories.} B.~Vallette in [19] introduces several versions
 of the notion of 2--monoidal category. This is a class of categories that could replace
 $(Vect_k,\otimes )$ in the following sense: operads with components in such
 a category could inherit black products, white products, and possibliy
 duality with the help of componentwise constructions, similar to the discussed ones.
 \smallskip

 Briefly, 2--monidal category  is a category with two monoidal
 structures (possibly non--commutative) say, $\otimes$ and $\boxtimes$, additionally endowed with a natural
 transformation called   {\it ``interchange law"}:
$$
 \varphi_{X,X^{\prime}, Y,Y^{\prime}}: \ (X\otimes X^{\prime})\boxtimes (Y\otimes Y^{\prime}) \to
 (X\boxtimes Y)\otimes (X^{\prime} \boxtimes Y^{\prime}) 
 $$
 satisfying certain compatibility conditions: see  [19], Proposition 2.
\smallskip

This interchange law replaces morphisms $S_{23}$ in (2.2) and similar situations.

%\smallskip

%It would be interesting to extend Vallette's formalism and to obtain new
%examples of Grothendieck--Verdier duality in  an operadic setting.

 \medskip

  \bigskip

 \centerline{\bf 3.  $F$--manifolds  and Dubrovin's almost duality}
 
 \medskip
 
 {\bf 3.1. Conventions and notations.} The basic definition of $F$--manifolds below
 works in each of the categories of 
manifolds $M$: $C^{\infty}$, analytic, or  formal,
eventually with odd (anticommuting 
coordinates). We denote the ground field
$k$, usually it is $\bold{C}$ or $\bold{R}$.

\smallskip

We denote the structure sheaf of $M$ by $\Cal{O}_M$ and  tangent sheaf by $\Cal{T}_M$.
The tangent sheaf is a locally free $\Cal{O}_M$--module;
its (super)rank is the (super)dimension of $M$.

\smallskip

 Let now $A$ be a linear $k$--(super)space with a $k$--bilinear commutative
multiplication and a $k$--bilinear Lie bracket.

\smallskip

We call its {\it Poisson tensor} the trilinear map for $a,b,c \in A$
$$
A^{\otimes 3}\ni a\otimes  b \otimes c \mapsto P_a(b,c):=[a,bc]-[a,b]c-(-1)^{ab}b[a,c].
$$
For $M$ as above, $\Cal{O}_M$ has a natural commutative multiplication, whereas  $\Cal{T}_M$
has a natural Lie structure. 
\smallskip
Poisson structure on $M$ involves introducing an extra Lie structure upon  $\Cal{O}_M$, whereas
$F$--structure involves  introducing an extra multiplication $\bullet$  upon $\Cal{T}_M$, satisfying the following
axiom.

\medskip

{\bf 3.2. Definition ([13]).} {\it A structure of $F$--manifold on $M$ is given by 
an $\Cal{O}_M$--bilinear associative commutative multiplication $\bullet$
on $\Cal{T}_M$ satisfying the so called $F$--identity:
$$
P_{X\bullet Y} =X\bullet P_Y+Y\bullet P_X.
\eqno(3.1)
$$
}
For brevity, we omitted here (obvious) signs relevant for the case of supermanifolds
and below will focus on the pure even case.
\smallskip

The geometric meaning of (3.1) was clarified in [14].
Namely, for any manifold $M$, consider
the sheaf of those functions on the 
cotangent manifold $T^*M$  which are polynomial
along the fibres of projection $T^*M\to M$. They constitute
the relative symmetric algebra $Symm_{\Cal{O}_M}(\Cal{T}_M)$

\smallskip

It is a sheaf of $\Cal{O}_M$--algebras, multiplication in which we denote $\bold{\cdot}$ \,.

\smallskip

Consider now a triple $(M,\bullet, e)$ where $\bullet$ is a commutative associative $\Cal{O}_M$--bilinear  multiplication
on $\Cal{T}_M$ with identity $e$.

\smallskip

There is an obvious homomorphism of $\Cal{O}_M$--algebras 
$$
(Symm_{\Cal{O}_M}(\Cal{T}_M), \cdot )\to (\Cal{T}_M, \bullet )
\eqno(3.2)
$$
\smallskip
{\bf 3.3. Theorem [14].}  {\it The multiplication $\bullet$ satisfies the $F$--identity (3.1)
iff the kernel  of (3.2) is stable with respect to 
the canonical Poisson brackets on $T^*M$.
\smallskip
In other words, $F$--identity is equivalent to the fact that the spectral cover
$$
\widetilde{M} := Spec_{\Cal{O}_M}(\Cal{T}_M, \bullet ) \to M
$$
(considered as 
a closed relative subscheme of the cotangent bundle of $M$) is coisotropic manifold of
 dimension $\r{dim}\, M$.}

\medskip

Notice that the spectral cover is not necessarily
a manifold. Its structure 
sheaf may have zero divisors and nilpotents ([14]).

\smallskip

However, it is a manifold, if the $F$--manifold $M$ is semisimple, which means that
locally $(\Cal{T}_M, \bullet )$ is isomorphic to a direct sum of $\r{dim}\, M$ copies of $\Cal{O}_M$.
 An embedded submanifold $N\subset T^*M$ is
the spectral cover of some semisimple $F$--structure
iff $N$ is Lagrangian.

\medskip

{\it Remark.} In the context of this note, it would be natural to define 
a general $F$--algebra as a linear space over a field $k$
endowed with a commutative algebra product $\bullet$ and  a Lie algebra bracket $[\, ,]$
that together satisfy the $F$--identity (3.1), and then to study the operad, classifying such algebras,
say $\O\Cal{F}$. This operad, however, does not fit in the context of 
quadratic dualities: the identity (3.1) written with additional arguments omitted in (3.1)
$$
[X\bullet Y,Z\bullet  W]-[X\bullet  Y,Z]\bullet  W-
Z\bullet  [X\bullet  Y,W]
$$
$$
-X\bullet  [Y,Z\bullet  W] + X\bullet  [Y,Z]\bullet  W+
X\bullet  Z\bullet  [Y,W]
$$
$$
-Y\bullet  [X,Z\bullet  W] + Y\bullet  [X,Z] \bullet W +
Y\bullet  Z\bullet  [X,W] =0
$$
operadically is a {\it  cubic relation of arity four.}

\medskip

{\bf 3.4. $F$--manifolds and mirror symmetry.}  One of the incarnations
of mirror symmetry involves isomorphisms of Frobenius manifolds
coming, say, from quantum cohomology, with Saito's Frobenius structures
on the germs of deformations of singularities.

\smallskip

By weakening Frobenius structure to $F$--structure, one can establish
a simple and beautiful class of examples of mirror symmetry.
Namely, we have  ([12], Theorems 5.3 and 5.6):

\medskip

{\bf 3.5. Theorem.} {\it  (i) The spectral cover space $\widetilde{M}$
of the canonical $F$--structure on the germ $M$ of the 
unfolding space
of an isolated hypersurface singularity is smooth.

\medskip

(ii) Conversely, let $M$ be an irreducible germ of a generically
semisimple $F$--manifold with 
the smooth spectral cover $\widetilde{M}.$
Then it is (isomorphic to) the germ of the unfolding space
of an 
isolated hypersurface singularity. Moreover, any isomorphism
of germs of such unfolding 
spaces compatible with their
$F$--structure comes from a  stable right equivalence of the
germs 
of the respective singularities.}

\smallskip

 Recall that the stable right equivalence is generated by adding sums
of squares of coordinates 
and making invertible local analytic
coordinate changes.

\smallskip

In view of this result, it would be important to understand
the following

\smallskip

{\bf Problem.} {\it Characterise those varieties $V$ for which the (genus zero)
 quantum cohomology Frobenius 
 spaces $H^*_{quant}(V)$ have
smooth spectral covers.}

\smallskip

Theorem 3.5 produces for such manifolds a weak version
of Landau--Ginzburg model, and 
thus gives a partial
solution of the mirror problem for them.

\medskip

{\bf 3.6. Dubrovin's duality for $F$--manifolds.} In [17], the following
version of  Dubrovin's almost duality ([10], see also [9]) was introduced:

\smallskip

{\bf Definition.} {\it An (even) vector field $\varepsilon$ on
an $F$--manifold with identity $(M,\bullet , e)$ is called
an eventual identity,
if $\varepsilon$ is $\bullet$--invertible, and moreover, the multiplication 
on vector fields
$$
X\circ Y:= X\bullet Y\bullet  \varepsilon^{-1}
$$
defines a new $F$--manifold structure with identity $(M, \circ , \varepsilon )$.}

\smallskip

There is a clear analogy between pairs of objects $(\bold{1}, K)$ considered in our discussion of 
GV--categories, and pairs of vector fields $(e,\varepsilon )$ on $M$, although
I do not know a natural definition of monodical category whose objects would be vector fields
and monoidal structure $\bullet$.
Perhaps, categorifications introduced in [18] and [8] might be
enlightening.

\smallskip

This analogy can be somewhat extended, using the following results of David and Strachan:

\medskip

{\bf 3.7. Theorem ([6]).} {\it (i) The field $\varepsilon$ is an eventual identity iff
for all $X,Y\in \Cal{T}_M$ we have
$$
P_{\varepsilon}(X,Y)= [e,\varepsilon ]\bullet X\bullet Y .
$$

(ii) On any $F$--manifold $(M,\bullet, e)$,
eventual identities form a group 
with respect to $\bullet$.
\smallskip
Moreover, if $\varepsilon_i$, $i=1,2$, are such eventual identities for  $(M,\bullet, e)$
that $[\varepsilon_1, \varepsilon_2]$ is invertible, 
then this commutator
is an eventual identity as well.}

\medskip

One may compare this result with Proposition 2.3 in [4] characterising the full
subcategory of dualizing objects (our $\varepsilon$'s) in a $GV$--category.

\bigskip

{\bf 3.8. Example 1.} For any eventual identity $\varepsilon$ on an $F$--manifold $(M,\bullet , e)$
and for any $m,n\in \bold{Z}$ we have
$$
[\varepsilon ^{\bullet n}, \varepsilon ^{\bullet m}] = (m-n)\varepsilon ^{\bullet (m+ n-1)}\bullet [e,\varepsilon ] .
$$

\bigskip

{\bf 3.9. Example 2.}  Let $(M,\bullet, e)$ be pure even and semisimple, with canonical
coordinates 
$(u^i)$, $i=1,\dots ,n$,  
so that  $\partial_i\bullet\partial_j=\delta_{ij}$,
where $\partial_i:=\partial/\partial u^i$.

\smallskip

Then eventual identities are precisely vector fields of the form
$$
\varepsilon =  f_1(u_1)\partial_1 +\dots  +f_n(u_n)\partial_n
$$
where $f_i$ are invertible functions of one variable. 

\smallskip

This statement can be compared with the discussion of  idempotents
in the $GV$--context in 3.4--3.7 of  [4].

\bigskip
\centerline{\bf References }

\medskip

[1]   M.~Barr. {\it *--autonomous categories.} Lecture Notes in Math., 752, 
Springer Verlag, 1979.

\smallskip

[2]    M.~Barr. {\it *--autonomous categories, revisited.} J.~Pure Appl. Algebra, 111(1996), 1--20.

%\smallskip

%[BaWe85] M.~Barr, C.~Wells. {\it Toposes, Triples and Theories.} Springer, 1985.

\smallskip
 
[3]   C.~Berger,  M.~Dubois--Violette, M.~Wambst. {\it Homogeneous
 algebras.} J. Algebra 261 (2003), no.~1, 172--185.

\smallskip

[4]  M.~Boyarchenko, V.~Drinfeld. {\it A duality formalism in the spirit 
 of Grothendieck and Verdier.} Quantum Topology, 4 (2013), 447--489.

 \smallskip
 
[5]     D.~Borisov, Yu.~Manin. {\it Generalized operads and their inner cohomomorphisms.}
 Progr. Math., vol 265, Birkh\"auser, Basel, 247--308.
 
% \smallskip
 
%[DavHe14]   L.~David, C.~Hertling. {\it Regular $F$--manifolds: initial conditions and Frobenius metrics.}
%arXiv:1411.4553v3

\smallskip
[6]   L.~David, I.~A.~B.~Strachan. {\it Dubrovin's duality for $F$--manifolds
with eventual identities.} Advances in Math. vol.~206 (2011), 4031--4060.
\smallskip
 
 [7]   B.~Day, R.~Street. {\it Quantum categories, star autonomy, and quantum groupoids.}
Galois theory, Hopf algebras, and semiabelian categories, 187Ð225, Fields Inst. 
Commun., 43, Amer. Math. Soc., Providence, RI, 2004. 

\smallskip

[8]   V.~Dotsenko, S.~Shadrin, B.~Vallette. {\it De Rham cohomology and
homotopy Frobenius manifolds.} Journ. Eur. Math. Soc.,vol. 17, no. 3 (2015). arXiv:1203.5077.

\smallskip

[9]   B.~Dubrovin. {\it Geometry of 2D topological field theory.}
In: Springer Lecture Notes in Math. 1620 (1996), 120--348.

\smallskip

[10]   B.~Dubrovin. {\it On almost duality for Frobenius manifolds.}
Geometry, topology and math phys., Amer. Math. Soc. Translations, Ser. 1, 212 (2004), 75--132.
arXiv: DG/0307374

 \smallskip
 
 [11]    V.~Ginzburg, M.~Kapranov. {\it Koszul duality for operads.}
 Duke Math. J. 76 (1994), no. 1, 203--272. {\it Erratum:}
 Duke Math. J. 80 (1995), no. 1, 293.
 
 \smallskip

[12]  C.~Hertling. {\it Frobenius manifolds and moduli spaces for singularities.}
Cambridge University Press, 2002.

\smallskip

[13]   C.~Hertling, Yu.~Manin. {\it Weak Frobenius manifolds.}
Int. Math. Res. Notices, 6 (1999), 277--286. 
arXiv:math.QA/9810132

\smallskip

[14]   C.~Hertling, Yu.~Manin, C.~Teleman. {\it An update on semisimple quantum cohomology
and $F$--manifolds.} Proc.~Steklov  Math.~Inst., vol.~264 (2009), 62--69. arXiv:math.AG/0803.2769

%\smallskip

%[LoMa00] A.~Losev, Yu.~Manin. {\it New moduli spaces of pointed curves and pencils of flat connections.} Fulton's %Festschrift,
%Michigan Journ. of Math., 48 (2000), 443--472. arXiv:math.AG/0001003

\smallskip
 
 [15]  J.-L.~Loday, B.~Vallette. {\it Algebraic operads.} Springer, 2012.

\smallskip
[16]  Yu.~Manin.  {\it Quantum groups and non--commutative geometry.}
CRM, Montr\'eal, 1988.  

\smallskip

[17]   Yu.~Manin. {\it F--manifolds with flat structure and
Dubrovin's duality.} Advances in Math. vol. 198:1 (2005), 5--26.  arXiv:math.DG/0402451

\smallskip

[18]  S.~Merkulov. {\it Operads, deformation theory and $F$--manifolds.} In: Frobenius Manifolds. Quantum
Cohomology and Singularities. Eds. C.~Hertling, M.~Marcolli. Aspects of Math. Vol.~ E36 (2004), 213--251.
arXiv:math/0210478

\smallskip

[19]  B.~Vallette. {\it Manin products, Koszul duality, Loday algebras
and Deligne conjecture.} J. Reine Angew. Math. 620 (2008), 105--164.
arXiv:math/0609002

   \enddocument